\documentclass[11pt]{article}
\usepackage[a4paper]{anysize}\marginsize{2.5cm}{2.5cm}{1.3cm}{1.2cm}
\pdfpagewidth=\paperwidth \pdfpageheight=\paperheight
\usepackage{pgf,tikz,hyperref}
\usetikzlibrary{arrows}
\usepackage{amsfonts,amssymb,amsthm,amsmath,eucal}
\usepackage{pdflscape}

\newcommand\numberthis{\addtocounter{equation}{1}\tag{\theequation}}

\theoremstyle{plain}
\newtheorem{thm}{Theorem}
\newtheorem{theorem}[thm]{Theorem}

\newtheorem{lemma}[thm]{Lemma}
\newtheorem{proposition}[thm]{Proposition}

\theoremstyle{definition}
\newtheorem{definition}[thm]{Definition}
\newtheorem{remark}[thm]{Remark}

\newtheorem{thevarthm}[thm]{\varthmname}

\newenvironment{varthm*}[1]{\trivlist\item[]{\bf #1.}\it}{\endtrivlist}

\newcommand\eps{\varepsilon}
\newcommand\be{\begin{eqnarray*}}
\newcommand\ee{\end{eqnarray*}}

\newcommand\C{\mathbb C}

\newcommand\cali{\mathcal I}

\renewcommand\P{\mathbb P}

\newcommand\calo{{\mathcal O}}
\newcommand\newop[2]{\def#1{\mathop{\rm #2}\nolimits}}
\newop\edim{edim}
\newop\Zeroes{Zeroes}
\newop\Jac{Jac}

\newcommand\keywords[1]{{\renewcommand\thefootnote{}\footnotetext{\textit{Keywords:} #1.}}}
\newcommand\subclass[1]{{\renewcommand\thefootnote{}\footnotetext{\textit{Mathematics Subject Classification (2010):} #1.}}}

\def\endproof{\hspace*{\fill}\endproofsymbol\endtrivlist}

\def\endproofsymbol{\frame{\rule[0pt]{0pt}{6pt}\rule[0pt]{6pt}{0pt}}}

\begin{document}

\author{Justyna~Szpond}
\title{Unexpected hypersurfaces with multiple fat points}
\date{\today}
\maketitle
\thispagestyle{empty}

\begin{abstract}
   Starting with the ground-breaking work of Cook II, Harbourne, Migliore
   and Nagel, there has been a lot of interest in unexpected hypersurfaces.
   In the last couple of months a considerable number of new examples
   and new phenomena has been observed and reported on.
   All examples studied so far had just one fat point. In this note
   we introduce a new series of examples, which establishes for the
   first time the existence of unexpected hypersurfaces with multiple
   fat points. The key underlying idea is to study Fermat-type
   configurations of points in projective spaces.
\subclass{MSC 14C20 \and MSC 14J26 \and MSC 14N20 \and MSC 13A15 \and MSC 13F20}
\end{abstract}


\section{Introduction}
   In the ground-breaking paper \cite{CHMN} the authors introduced the notion of unexpected curves.
\begin{definition}\label{def: unexpected curve}
   We say that a reduced finite set of points $Z\subset \P^2$ \emph{admits an unexpected curve} of degree
   $j+1$ if there is an integer $j>0$ such that, for a general point $P$, $jP$ fails to impose the expected
   number of conditions on the linear system of curves of degree $j+1$ containing $Z$. That is, $Z$ admits
   an unexpected curve of degree $j+1$ if
$$
   h^0(\cali_{Z+jP}(j+1))>\max\left\{ h^0(\cali_Z(j+1))-\binom{j+1}{2},0\right\}.
$$
\end{definition}
   This notion has been generalized to hypersurfaces in projective spaces of arbitrary dimension
   in the subsequent paper \cite{HMNT}.
\begin{definition}\label{def: unexpected hypersurface}
   We say that a reduced set of points $Z\subset\P^N$ \emph{admits an unexpected hypersurface} of degree $d$
   if there exists a sequence of integers $m_1,\ldots,m_s$ such that for general points $P_1,\ldots,P_s$
   the zero-dimensional subscheme $m_1P_1+\ldots +m_sP_s$ fails to impose independent conditions on forms
   of degree $d$ vanishing along $Z$.
\end{definition}
   A well known example of this kind is provided by $7$ general double points in $\P^4$.
   There is a twisted quartic curve passing through these $7$ points and its secant variety
   is a threefold of degree $3$ vanishing to order $2$ along the quartic curve, in particular
   in the $7$ general points. This threefold is unexpected, see \cite{AleHir95}.

   This example has two deficiencies. First, the singular points are not isolated. Second, the set $Z$ is empty.

   The first example of an unexpected surface $Q_R$ in $\P^3$ admitted by a non-empty set $Z$
   has been discovered and described in \cite[Theorem 1]{BMSS}.
   In that example there is an isolated point $R$ of multiplicity $3$ which is a general point.
   The surface $Q_R$ has $4$ additional double points, whose coordinates depend on $R$, and no other singularities.
   In the present note we show the following result.
\begin{varthm*}{Main Theorem}
   There exists a non-empty set $Z$ in $P^N$ with $N=2k+1$ for $k\geq 1$ which admits an
   unexpected hypersurface for multiplicities $3,\underbrace{2,\ldots,2}_{k-1}$.
   This hypersurface has only isolated singularities.
\end{varthm*}
   For $N=3$ we have a computer-free proof. Our construction in higher dimensional projective
   spaces builds upon the existence of an unexpected surface of degree $4$ in $\P^3$ with
   a general point of multiplicity $3$. Given a general point $R$ in $\P^5$ we project it to
   a coordinate $\P^3$ (i.e. codimension $2$ flat $\Pi$ defined by the vanishing of two coordinates)
   obtaining the point $R_{\Pi}$.
   Then the $N=3$ case gives an unexpected surface in $\Pi$ which vanishes to order $3$ at $R_{\Pi}$.
   Taking a cone in $\P^5$ with vertex in the coordinate line spanned by the coordinate points
   corresponding to the two coordinates which defined $\Pi$ we get a $4$-fold of degree $4$ in $\P^5$ which
   in particular vanishes to order $3$ at $R$. Of course this $4$-fold has non-isolated singularities.
   However, taking the linear system of quartics generated by all such cones in $\P^5$
   we are able to identify its $6$ generators. Given another general point $P$ in $\P^5$ we find
   an explicit formula for coefficients
   in front of the $6$ generators (depending on $P$) so that the resulting $4$-fold has only isolated
   singularities and in particular has multiplicity $3$ at $R$ and multiplicity $2$ at $P$.
   Once the equation is in place, all our claims can be (in principle) checked by direct
   computation. We provide a Singular script \cite{jusF} to ease their verification. The construction
   in higher dimensions is performed along the same lines. So far it justified only by computer
   computations with random points rather than arbitrary as in the case of Theorem \ref{thm:P5 with mult 3 and 2}.

   The main idea here is to explore Fermat-type configurations of points in projective spaces,
   see \cite{Szp19c} for a recent survey on Fermat-type hyperplane arrangements and their
   applications in algebra and geometry.

\section{Fermat-type point configurations}\label{sec:Fermat}
   By the way of warm-up we begin with points in $\P^2$. In Theorem \ref{thm: unexpected curves mult 4}
   we show a peculiar family of unexpected curves, which degree grows while the multiplicity in the
   general point remains fixed, equal $4$. Passing then to the general case of $\P^N$ we show identify
   generators of the ideal of Fermat-type configuration of points. This is crucial for the subsequent part.
\subsection{Fermat point configurations in $\P^2$}
   Let $n\geq 1$ be a positive integer and let $I_{2,n}$ be the complete intersection ideal in $\C[x,y,z]$
   viewed as the graded ring of $\P^2$ generated by
   \be
      x^n-y^n,\
      y^n-z^n.
   \ee
   Let $F_{2,n}$ be the ideal
   $$F_{2,n}=I_{2,n}\cap (x,y)\cap (x,z)\cap (y,z).$$
   The support of $I_{2,n}$ is the set of $n^2$ points
   \be
      (1:\eps^\alpha:\eps^\beta)
   \ee
   where $\eps$ is a primitive root of unity or order $n$ and
   $1\le\alpha,\beta\le n$. The support of $F_{2,n}$ is the union
   of the support of $I_{2,n}$ and the $3$ coordinate points
   $(1:0:0)$, $(0:1:0)$ and $(0:0:1)$.

   The ideals $F_{2,n}$ have appeared recently in various guises, most notably in the context
   of the containment problem, see \cite{SzeSzp17} for an introduction to this circle of ideas
   and \cite{MalSzp18Moieciu}, \cite{MalSzp18} for specific applications.
   Their algebraic properties have been studied in depth by Nagel and Seceleanu \cite{NagSec16}.

\begin{figure}
\centering
\begin{tikzpicture}[line cap=round,line join=round,x=1.0cm,y=1.0cm,scale=0.7]
\clip(1,-2) rectangle (9,4.5); \draw (1.5,-1.5)-- (5.25,2.25); \draw
(2.5,2.)-- (5.,2.); \draw (2.5,1.)-- (6.75,1.); \draw (2.5,0.)--
(5.,0.); \draw (3.,2.)-- (3.,-0.5); \draw (4.,3.75)-- (4.,-0.5);
\draw (5.,2.)-- (5.,-0.5); \draw  plot [smooth, tension=2]
coordinates {(3,1)  (2,-1) (5,0)}; \draw  plot [smooth, tension=1]
coordinates {(3,2) (1.3, 0.2)  (2,-1) (4,0)}; \draw (3.,1.)--
(4.25,2.25); \draw (4.,0) -- (5.25,1.25); \draw (3.,1.)--
(4.25,2.25); \draw (3.,2.)-- (3.25,2.25); \draw (5.,0.)--
(5.25,0.25); \draw [fill=black] (3,2) circle (1.5pt); \draw
[fill=black] (4,2) circle (1.5pt); \draw [fill=black] (5,2) circle
(1.5pt); \draw [fill=black] (3,1) circle (1.5pt); \draw [fill=black]
(4,1) circle (1.5pt); \draw [fill=black] (5,1) circle (1.5pt); \draw
[fill=black] (3,0) circle (1.5pt); \draw [fill=black] (4,0) circle
(1.5pt); \draw [fill=black] (5,0) circle (1.5pt); \draw [fill=black]
(2,-1) circle (1.5pt); \draw [fill=black] (6.5,1) circle (1.5pt);
\draw [fill=black] (4,3.5) circle (1.5pt); \draw  plot [smooth,
tension=2] coordinates {(3.,2)  (3.38,2.75) (4,3.5) }; \draw
(4,3.5)-- (4.25,3.75); \draw  plot [smooth, tension=2] coordinates
{(5,2)  (4.62,2.75) (4,3.5) }; \draw (4,3.5)-- (3.75,3.75); \draw
plot [smooth, tension=2] coordinates {(5.,2)  (5.62,1.7) (6.5,1) };
\draw (6.5,1)-- (6.75,1.25); \draw  plot [smooth, tension=2]
coordinates {(5.,0)  (5.62,0.3) (6.5,1) }; \draw (6.5,1)--
(6.75,0.75);
\end{tikzpicture}
\caption{Dual Hesse arrangement}
\label{fig: dual Hesse}
\end{figure}
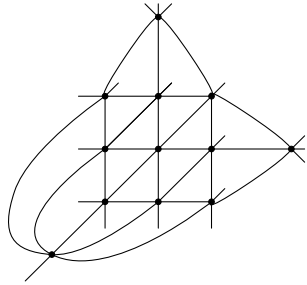

   For $n=3$ the set of points defined by $F_{2,3}$ is the set of singular points
   (i.e. points where two or more arrangement lines meet) of the dual Hesse arrangement
   of $9$ lines with $12$ triple points at the zeroes of $F_{2,3}$. This arrangement
   is depicted in Figure \ref{fig: dual Hesse}. The lines are indicated by curved segments
   as it is not possible to embed this configuration in the real projective plane due to
   the celebrated Sylvester-Gallai Theorem, see \cite{BorweinMoser1990}.
   It has been proved already in \cite[Lemma 2.1]{DST13} that the ideal $F_{2,3}$ is generated in degree
   $4$. More precisely we have
\begin{lemma}\label{lem:generators in P2}
   The ideal $F_{2,3}$ is generated by the binomials
   $$x(y^3-z^3),\;
     y(z^3-x^3),\;
     z(x^3-y^3).$$
\end{lemma}
   In fact exactly the same proof works for any $n\geq 3$, so that we have
\begin{lemma}\label{lem: generators in P2 n}
   The ideal $F_{2,n}$ is generated by the binomials
   $$x(y^n-z^n),\;
     y(z^n-x^n),\;
     z(x^n-y^n).$$
\end{lemma}
   This ideal is supported on $n^2$ points which form a complete intersection and the $3$ coordinate points.
   It turns out that this set of points admits an unexpected curve (Definition \ref{def: unexpected curve})
   for any $n\geq 3$.
\begin{theorem}[Unexpected curves with a point of multiplicity $4$]\label{thm: unexpected curves mult 4}
   Let $Z$ be the configuration of points in $\P^2$ defined by $F_{2,n}$
   for $n\geq 3$.
   Let $R=(a:b:c)$ be a general point in $\P^2$.
   With
   $u=\binom{n}{2}-1, \quad v=\binom{n-1}{2}, \quad w=\binom{n+1}{2}$,
   the polynomial
   \begin{align*}
   Q_P(x:y:z)=&-cxy((ub^n+vc^n)(z^n-x^n)+(ua^n+vc^n)(y^n-z^n))\\
        &-bxz((ua^n+vb^n)(y^n-z^n)+(uc^n+vb^n)(x^n-y^n))\\
        &-ayz((ub^n+va^n)(z^n-x^n)+(uc^n+va^n)(x^n-y^n))\\
        &+wa^{n-1}bcx^2(y^n-z^n)+wab^{n-1}cy^2(z^n-x^n)\\
        &+wabc^{n-1}z^2(x^n-y^n) \numberthis\label{eq: mult 4}
   \end{align*}
   \begin{itemize}
   \item vanishes at all points of $Z$,
   \item vanishes to order $4$ at $R$,
   \item defines an unexpected curve of degree $n+2$ for $Z$ with respect to $R$.
   \end{itemize}
\end{theorem}
\proof
   This claim is easy verify by direct calculations, which we omit.
\endproof
   Theorem \ref{thm: unexpected curves mult 4} is of interest as it exhibits a first family of unexpected curves,
   where the degree of curves grows but the multiplicity at the general point remains constant.
\subsection{Fermat-type configurations of points in $\P^N$}
   For a positive integer $n$ let $I_{N,n}$ be the complete intersection ideal defined
   by the binomials
   \begin{equation}\label{eq:def of I_N,n}
      x_0^n-x_1^n,\; x_1^n-x_2^n,\;\ldots,\; x_{N-1}^n-x_N^n.
   \end{equation}
   The set $Z_{N,n}$ of zeroes of $I_{N,n}$ is the set of $N^n$ points of the form
   $$(1:\eps^{\alpha_1}:\eps^{\alpha_2}:\ldots:\eps^{\alpha_N}),$$
   where $\eps$ is a primitive root of unity of order $n$ and
   $$1\leq\alpha_1,\alpha_2,\ldots,\alpha_N\leq n.$$
   Adding to this set of points the coordinate points of $\P^N$, we obtain
   the set $W_{N,n}$ defined by the ideal $F_{N,n}$.

   In the sequel
   we will deal with the index $n=3$, so we drop the degree from the notation.
   We also use the convention that the indices are understood modulo $N+1$.
   For example $x_{N+2}=x_1$. Moreover it is convenient to introduce the
   following notation
   $$[a,b]=a^3-b^3,$$
   which we use for numbers and for monomials. Note that the symbol $[a,b]$
   is anti-symmetric and satisfies a Jacobi-type identity
   \begin{equation}\label{eq:Jacobi type}
      a^3[b,c]+b^3[c,a]+c^3[a,b]=0\;\;\mbox{ for all }\; a,b,c.
   \end{equation}

   Lemma \ref{lem:generators in P2} generalizes in the following way.
\begin{lemma}\label{lem:generators in PN}
   The ideal $F_N$ is generated by the binomials of the form
   \begin{equation}\label{eq:generators in PN}
      x_i[x_{i+1},x_j],
   \end{equation}
   where $i\in\left\{0,1,\ldots, N\right\}$ and
   $j\in \left\{0,1,\ldots, N\right\}\setminus\left\{i,i+1\right\}$.
\end{lemma}
\proof
   Let $J$ be the ideal generated by binomials in \eqref{eq:generators in PN}.
   It is clear that elements in $J$ vanish in all points of $W_N$,
   so there is $J\subset F_N$.

   It is also easy to see that $J$ contains
   all binomials of the form $x_i[x_j,x_k]$, where the indices are mutually distinct. Indeed, we have
   $$x_i[x_j,x_k]=x_i[x_{i+1},x_k]-x_i[x_{i+1},x_j].$$
   It is also clear that the ideal $I_{N}$ defined in \eqref{eq:def of I_N,n} is generated also by binomials
   $$[x_0,x_1],\;[x_0,x_2],\;\ldots,\;[x_0,x_N],$$
   which are slightly more convenient to work with in this proof.

   Let $f\in F_N$ be an arbitrary polynomial. We want to show the containment $f\in J$ or equivalently $f=0 \mod J$.
   Since $f$ vanishes in particular in all points of $Z_N$, there are homogeneous polynomials $g_1,\ldots,g_N$ such that
   \begin{equation}\label{eq:f=}
      f=\sum_{i=1}^{N} g_i[x_0,x_i].
   \end{equation}
   from now on, we work$\mod J$. Since for $j\neq 0,i$ the binomials $x_j[x_0,x_i]$ are in $J$,
   we may assume that each polynomial $g_i$ depends only on $x_0$ and $x_{i}$. Moreover, for a fixed
   $j\neq 0, i$ we have
   $$x_0x_i[x_0,x_i]=x_i\cdot x_0[x_j,x_i]-x_0\cdot x_i[x_j,x_0]\in J.$$
   Thus it must be
   $$g_i=a_ix_0^d+b_ix_i^d \mod J$$
   for some $d\geq 0$ and some scalars $a_1,\ldots, a_N,b_1,\ldots, b_N\in\C$.
   Evaluating \eqref{eq:f=} in the coordinate points we obtain
   \begin{equation}\label{eq:a_i b_i}
      a_1+a_2+\ldots+a_N=0\mbox{ and } b_1=b_2=\ldots=b_N=0.
   \end{equation}
   Hence \eqref{eq:f=} reduces to
   \begin{equation}\label{eq:f= 2}
      f=a_1x_0^d[x_0,x_1]+a_2x_0^d[x_0,x_2]+\ldots+a_Nx_0^d[x_0,x_N] \mod J.
   \end{equation}
   If $d=0$, then evaluating \eqref{eq:f= 2} in the coordinate points, we obtain
   $$a_1=a_2=\ldots=a_N=0$$
   and we are done.

   If $d>0$, then, using first part of \eqref{eq:a_i b_i}, we can rewrite \eqref{eq:f= 2} as
   $$ f=a_1x_0^d [x_N,x_1]+a_2x_0^d[x_N,x_2]+\ldots+a_{N-1}x_0^d[x_N,x_{N-1}],$$
   which is clearly an element of $J$ and we are done again.
\endproof

\section{Unexpected quartic surfaces in $\P^3$ with a triple point}
   The story of unexpected hypersurfaces begins with unexpected curves
   discovered by Di Gennnro, Ilardi and Vall\`es in \cite{DIV14}
   and studied systematically by Cook II, Harbourne, Migliore and Nagel in \cite{CHMN}.
   In \cite{BMSS} the first example of a higher dimensional unexpected hypersurface (a surface
   in $\P^3$) has been described by Bauer, Malara, Szemberg and the present author.
   Shortly after that Harbourne, Migliore, Nagel and Teitler in \cite{HMNT}
   constructed examples of unexpected hypersurfaces in projective spaces
   of arbitrary dimension. All these examples have just one singular general point.

   Here we recall the construction of the unexpected quartic in $\P^3$.
   This example has not only motivated our construction of unexpected hypersurfaces
   with multiple fat points but it serves as a building fundament for this construction.

   According to Lemma \ref{lem:generators in PN}, in $\P^3$ the ideal $F_{3,3}$
   has $8$ generators:
   $$g_{0,2}=x_0[x_1,x_2],\;
     g_{0,3}=x_0[x_1,x_3],\;
     g_{1,3}=x_1[x_2,x_3],\;
     g_{1,0}=x_1[x_2,x_0],$$
   $$g_{2,0}=x_2[x_3,x_0],\;
     g_{2,1}=x_2[x_3,x_1],\;
     g_{3,1}=x_3[x_0,x_1],\;
     g_{3,2}=x_3[x_0,x_2].$$
   It is convenient to introduce the following notation. Given mutually
   distinct numbers $i,j\in\left\{0,1,2,3\right\}$, we denote by $k$ the
   index $\left\{0,1,2,3\right\}\setminus\left\{i,i+1,j\right\}$. As usually,
   the indices are considered modulo $4$.
   With this notation we have (compare \cite[Theorem 6]{BMSS}).
\begin{theorem}\label{thm:quartic in P3}
   For a general point $R=(a_0:a_1:a_2:a_3)$,
   the quartic
   $$Q_R(x_0:x_1:x_2:x_3)=\sum\limits_{i=0}^3\sum\limits_{j=i+2}^{i+3}(-1)^k a_i^2 [a_{i+1},a_k]\cdot g_{i,j}$$
   has a triple point at $R$ and satisfies $Q_R\in F_{3,3}$.
\end{theorem}
\proof
   The second property is obvious. Checking that $Q_R$ vanishes to order $3$ at $R$
   is a straightforward computation. However for further reference we will present
   some computations which make heavy use of available symmetries.

   To begin with, we check that the first order derivatives of $Q_R$ vanish at $R$.
   By symmetry it is enough to check just one derivative. We compute the one with respect to $x_0$
   and obtain
   $$\frac{\partial}{\partial x_0}Q_R=\left([a_1a_2][x_1x_3]-[a_1a_3][x_1x_2]\right)a_0^2+
     3x_0^2\left(a_1^2x_1[a_2a_3]+a_2^2x_2[a_3a_1]+a_3^2x_3[a_1a_2]\right).$$
   Vanishing of the first bracket when evaluating at $R$ is straightforward. The vanishing of the second
   bracket follows from the Jacobi-type identity \eqref{eq:Jacobi type}.

   For the second order derivative with respect to $x_0$ we have
   $$\frac{\partial^2}{\partial x_0^2}Q_R=
     -6x_0\cdot\sum_{k=1}^{3}(-1)^ka_j^2a_l^2(a_lx_j-a_jx_l),$$
   where the indices $j,l$ are chosen in dependence on $k$ so that
   $j<l$ and $\left\{j,k,l\right\}=\left\{1,2,3\right\}$.

   For the mixed derivative of order $2$ we may again use the symmetry
   and look at just one other variable, here $x_1$. We have
   $$\frac{\partial^2}{\partial x_1 \partial x_0}Q_R=
     3[a_2a_3](a_1^2x_0^2-a_0^2x_1^2).$$

   Clearly, in both cases we obtain zero, when evaluating at $R$.

   Instead of computing derivatives one by one, we can also use the
   whole interpolation matrix at once. We present this matrix in Table \ref{tab: imat 3},
   since this is a baby case for our argument in the proof of Theorem \ref{thm:P5 with mult 3 and 2},
   see Table \ref{tab: imat 5}.
\begin{table}
\begin{equation*}
  \begin{pmatrix}
   0 &  0 &  -a_1 &  0 &  -a_2 &  0 &  a_3 &  a_3 \\
   a_1^2 &  a_1^2 &  -a_0^2 &  0 &  0 &  0 &  0 &  0 \\
   a_2^2 &  0 &  0 &  0 &  a_0^2 &  0 &  0 &  0 \\
   0 &  -a_3^2 &  0 &  0 &  0 &  0 &  a_0^2 &  a_0^2 \\
   a_0 &  a_0 &  0 &  0 &  0 &  -a_2 &  -a_3 &  0 \\
   0 &  0 &  a_2^2 &  a_2^2 &  0 &  -a_1^2 &  0 &  0 \\
   0 &  0 &  0 &  a_3^2 &  0 &  0 &  a_1^2 &  0 \\
   -a_0 &  0 &  a_1 &  a_1 &  0 &  0 &  0 &  -a_3 \\
   0 &  0 &  0 &  0 &  a_3^2 &  a_3^2 &  0 &  -a_2^2 \\
   0 &  -a_0 &  0 &  -a_1 &  a_2 &  a_2 &  0 &  0
  \end{pmatrix}
\end{equation*}
\caption{Interpolation matrix for a point of multiplicity $3$ in $\P^3$.}
\label{tab: imat 3}
\end{table}
   The columns of the table correspond to the generators of $F_{3,3}$, whereas
   the rows stand for order $2$ partial derivatives. Taking into account that the point $R$
   is general, we have divided out common factors appearing in some rows. This is the reason
   why not all rows have entries of the same degree. Nevertheless, the matrix as entire, remains
   of course homogeneous.
\endproof
\section{Unexpected quartic $4$--folds in $\P^5$ with a triple and double point}
   In this section we present in detail unexpected hypersurfaces in $\P^5$,
   which have $2$ general fat points, one of multiplicity $3$ and the other one
   of multiplicity $2$.

   According to Lemma \ref{lem:generators in PN} the ideal $F_{5,3}$ has $24$ generators
   and its zero locus consists of $3^5$ Fermat points and $6$ coordinate points, so that
   there are altogether $249$ points in $W_{5,3}$. The linear system of quartics in $\P^5$
   has (affine) dimension $126$, so the points in $W_{5,3}$ do not impose independent conditions
   on quartics. Of course there exists a subset of $102$ points in $W_{5,3}$ which impose
   independent conditions on quartics in $\P^5$ but to determine such a set is not relevant
   for our considerations here.

   Since vanishing to order $3$ imposes $21$ conditions, it is expected that the system
   of quartics in $\P^5$ vanishing along $W_{5,3}$ and having a triple point in a general
   point $R=(a_0:a_1:\ldots:a_5)$ has dimension $3$. It turns out however that this system
   has in fact higher dimension.
\begin{proposition}\label{prop:P5 with mult 3}
   Let $I(R)$ denote the ideal of the point $R$. For
   $$V_{5;3}=H^0(\P^5; \calo_{\P^5}(4)\otimes F_{5,3}\otimes I(R)^3)$$
   we have $\dim V_{5;3}=6$.
\end{proposition}
\proof
   We verify our claim computing the interpolation matrix at $R$. More precisely, we
   compute all partial derivatives of order $2$ (there are $21$ of them) of all $24$
   generators of $F_{5,3}$. The matrix we get has a lot of zeros and thus its rank
   is easy to determine. The rank is the number of \emph{independent} conditions
   imposed on quartics in $F_{5,3}$ by vanishing at $R$ to order $3$.
\begin{landscape}
\begin{table}
\begin{equation*}
  \begin{pmatrix}
    0 & 0 & 0 & 0 & -a_1 & 0 & 0 & 0 & -a_2 & 0 & 0 & 0 & -a_3 & 0 & 0 & 0 & -a_4 & 0 & 0 & 0 & a_5 & a_5 & a_5 & a_5 \\
    a_1^2 & a_1^2 & a_1^2 & a_1^2 & -a_0^2 & 0 & 0 & 0 & 0 & 0 & 0 & 0 & 0 & 0 & 0 & 0 & 0 & 0 & 0 & 0 & 0 & 0 & 0 & 0 \\
    a_2^2 & 0 & 0 & 0 & 0 & 0 & 0 & 0 & a_0^2 & 0 & 0 & 0 & 0 & 0 & 0 & 0 & 0 & 0 & 0 & 0 & 0 & 0 & 0 & 0 \\
    0 &  a_3^2 &  0 &  0 &  0 &  0 &  0 &  0 &  0 &  0 &  0 &  0 &  a_0^2 &  0 &  0 &  0 &  0 &  0 &  0 &  0 &  0 &  0 &  0 &  0 \\
    0 & 0 & a_4^2 & 0 & 0 & 0 & 0 & 0 & 0 & 0 & 0 & 0 & 0 & 0 & 0 & 0 &  a_0^2 & 0 & 0 & 0 & 0 & 0 & 0 & 0 \\
    0 & 0 & 0 & -a_5^2 & 0 & 0 & 0 & 0 & 0 & 0 & 0 & 0 & 0 & 0 & 0 & 0 & 0 & 0 & 0 & 0 & a_0^2 &  a_0^2 &  a_0^2 & a_0^2 \\
    a_0 &  a_0 &  a_0 &  a_0 & 0 & 0 & 0 & 0 & 0 & -a_2 & 0 & 0 & 0 & -a_3 & 0 & 0 & 0 & -a_4 & 0 & 0 & -a_5 & 0 & 0 & 0 \\
    0 & 0 & 0 & 0 &  a_2^2 &  a_2^2 & a_2^2 &  a_2^2 & 0 & -a_1^2 & 0 & 0 & 0 & 0 & 0 & 0 & 0 & 0 & 0 & 0 & 0 & 0 & 0 & 0 \\
    0 & 0 & 0 & 0 & 0 &  a_3^2 & 0 & 0 & 0 & 0 & 0 & 0 & 0 &  a_1^2 & 0 & 0 & 0 & 0 & 0 & 0 & 0 & 0 & 0 & 0 \\
    0 &  0 &  0 &  0 &  0 &  0 &  a_4^2 &  0 &  0 &  0 &  0 &  0 &  0 &  0 &  0 &  0 &  0 &  a_1^2 &  0 &  0 &  0 &  0 &  0 & 0 \\
    0 &  0 &  0 &  0 &  0 &  0 &  0 &  a_5^2 &  0 &  0 &  0 &  0 &  0 &  0 &  0 &  0 &  0 &  0 &  0 &  0 &  a_1^2 &  0 &  0 & 0 \\
    -a_0 &  0 &  0 &  0 &  a_1 &  a_1 &  a_1 &  a_1 &  0 &  0 &  0 &  0 &  0 &  0 &  -a_3 &  0 &  0 &  0 &  -a_4 &  0 &  0 &  -a_5 &  0 & 0 \\
    0 &  0 &  0 &  0 &  0 &  0 &  0 &  0 &  a_3^2 &  a_3^2 &  a_3^2 &  a_3^2 &  0 &  0 &  -a_2^2 &  0 &  0 &  0 &  0 &  0 &  0 &  0 &  0 & 0 \\
    0 &  0 &  0 &  0 &  0 &  0 &  0 &  0 &  0 &  0 &  a_4^2 &  0 &  0 &  0 &  0 &  0 &  0 &  0 &  a_2^2 &  0 &  0 &  0 &  0 & 0 \\
    0 &  0 &  0 &  0 &  0 &  0 &  0 &  0 &  0 &  0 &  0 &  a_5^2 &  0 &  0 &  0 &  0 &  0 &  0 &  0 &  0 &  0 &  a_2^2 &  0 & 0 \\
    0 &  -a_0 &  0 &  0 &  0 &  -a_1 &  0 &  0 &  a_2 &  a_2 &  a_2 &  a_2 &  0 &  0 &  0 &  0 &  0 &  0 &  0 &  -a_4 &  0 &  0 &  -a_5 & 0 \\
    0 &  0 &  0 &  0 &  0 &  0 &  0 &  0 &  0 &  0 &  0 &  0 &  a_4^2 &  a_4^2 &  a_4^2 &  a_4^2 &  0 &  0 &  0 &  -a_3^2 &  0 &  0 &  0 & 0 \\
    0 &  0 &  0 &  0 &  0 &  0 &  0 &  0 &  0 &  0 &  0 &  0 &  0 &  0 &  0 &  a_5^2 &  0 &  0 &  0 &  0 &  0 &  0 &  a_3^2 & 0 \\
    0 &  0 &  -a_0 &  0 &  0 &  0 &  -a_1 &  0 &  0 &  0 &  -a_2 &  0 &  a_3 &  a_3 &  a_3 &  a_3 &  0 &  0 &  0 &  0 &  0 &  0 &  0 & -a_5 \\
    0 & 0 & 0 & 0 & 0 & 0 & 0 & 0 & 0 & 0 & 0 & 0 & 0 & 0 & 0 & 0 & a_5^2 & a_5^2 & a_5^2 & a_5^2 & 0  & 0 &  0 & -a_4^2 \\
    0 &  0 &  0 &  -a_0 &  0 &  0 &  0 &  -a_1 &  0 &  0 &  0 &  -a_2 &  0 &  0 &  0 &  -a_3 &  a_4 &  a_4 &  a_4 &  a_4 &  0 &  0 &  0 &  0
  \end{pmatrix}
\end{equation*}
\caption{Interpolation matrix for a point of multiplicity $3$ in $\P^5$.}
\label{tab: imat 5}
\end{table}
\endproof
\end{landscape}
   Since a double point in $\P^5$ imposes $6$ conditions on forms of arbitrary degree,
   we do not expect that there is a non-zero section in $V_{5;3}$ vanishing to order $2$
   at an additional general point. However such unexpected hypersurface does exist!

   We define first the following $6$ cones in $\P^5$. For indices $i,j\in\left\{0,\ldots,5\right\}$,
   with $i<j$, we denote by $s,t,u,v$ the remaining $4$ indices in the growing order. Then, for
   $i,j\in\left\{0,1,2,3\right\}$ with $i<j$ we set
   $$R_{i,j}=(a_s: a_t: a_u: a_v)$$ and
   $$J_{i,j}=Q_{R_{i,j}}(x_s: x_t: x_u: x_v),$$
   where $Q_R$ is taken from Theorem \ref{thm:quartic in P3}.
\begin{remark}
   It can be shown that $J_{0,1},\ldots,J_{2,3}$ generate $V_{5,3}$ but this is irrelevant
   for the further assertions and we omit the verification of this claim.
\end{remark}
\begin{theorem}\label{thm:P5 with mult 3 and 2}
   Let $P=(b_0:b_1:\ldots:b_5)$ be a general point in $\P^5$. Then there exists a unique
   quartic $Q_{R,P}\in V_{5;3}$ vanishing at
   \begin{itemize}
   \item all points of the Fermat-type configuration $W_5$,
   \item the point $R=(a_0:a_1:\ldots:a_5)$ to order $3$,
   \item the point $P$ to order $2$.
   \end{itemize}
\end{theorem}
\proof
   We are able to write down explicit equation of the quartic
   \begin{align*}
      Q_{R,P}(x_0:\ldots:x_5) &= J_{2,3}(P)\cdot J_{0,1}(x_0:\ldots:x_5)
                               - J_{1,3}(P)\cdot J_{0,2}(x_0:\ldots:x_5)\\
                               &+ J_{0,3}(P)\cdot J_{1,2}(x_0:\ldots:x_5)
                               + J_{1,2}(P)\cdot J_{0,3}(x_0:\ldots:x_5)\\
                               &- J_{0,2}(P)\cdot J_{1,3}(x_0:\ldots:x_5)
                               + J_{0,1}(P)\cdot J_{2,3}(x_0:\ldots:x_5).
   \end{align*}
   The sign in the front of a summand of $Q_{R,P}$ depends on whether there
   is a pair of consecutive numbers in the indices of involved $J$'s or not.
   If there isn't, i.e., for pairs $(0,2)$ and $(1,3)$ we get a minus.

   Note that the first two assertions of the Theorem follow straightforward
   from the way we defined the cones $J_{0,1},\ldots,J_{2,3}$. The last
   assertion can be checked by direct computations as in the proof of Theorem \ref{thm:quartic in P3}.

   The most tricky part of the proof was to find the coefficients in front of the cones.
\endproof
\section{Unexpected quartic $2k$-folds in $\P^{2k+1}$ with a triple point and $k-1$ double points}
   Before we conclude we collect a couple of facts towards the proof of our Main Theorem.
   We begin with an observation which parallels Proposition \ref{prop:P5 with mult 3}
\begin{proposition}\label{prop: PN with mult 3}
   Let $N=2k+1$ for some $k\geq 1$. Let $R=(a_0:\ldots:a_N)$ be a general point in $\P^N$ and let $I(R)$ be its saturated ideal.
   Let
   $$V_{N,3}=H^0(\P^N;\calo_{\P^N}(4)\otimes F_{N,3}\otimes I(R)^3).$$
   Then $\dim V_{N,3}=\binom{N-1}{2}$.
\end{proposition}
\begin{remark}
   The statement of Proposition \ref{prop: PN with mult 3} means that a general triple point always imposes exactly $3$
   conditions less on generators of $F_{N,3}$ than expected. We don't know how to verify this theoretically.
\end{remark}
   The final step in proving Main Theorem is the next statement.
\begin{proposition}
   Let $P_1,\ldots,P_{k-1}$ be general points in $\P^N$ with $N=2k+1$. Vanishing to order $2$ in these points
   imposes $2k^2-k-1$ conditions on elements in $V_{N,3}$.
\end{proposition}
\begin{remark}
   Note that vanishing to order $2$ at a general point is expected to impose $N+1=2k+2$ conditions.
   Thus $(k-1)$ points are expected to impose $2k^2-2$ conditions. Interestingly, computer experiments
   show that the first $(k-2)$ points do impose on $V_{N,3}$ the expected number of conditions, namely
   $2k^2-2k-4$. It is only the last point which imposes $k+3$ instead of $2k+2$ conditions
   on the system
   $$H^0(\P^N; \calo_{\P^N}(4)\otimes F_{N,3}\otimes I(R)^3\otimes I(P_1)^2\otimes\ldots\otimes I(P_{k-2})^2).$$
   We do not have a theoretical explanation of this phenomena.
\end{remark}




\bigskip
\bigskip
\small

   Justyna Szpond,
   Department of Mathematics, Pedagogical University of Cracow,
   Podchor\c a\.zych 2,
   PL-30-084 Krak\'ow, Poland.

\nopagebreak
   \textit{E-mail address:} \texttt{szpond@gmail.com}


\end{document}